\documentclass[fleqn]{zzz}
\usepackage{mathrsfs}
\usepackage{color}
\usepackage{tikz}
\usepackage{amsmath,amsthm,amsbsy,amssymb}
\usepackage{graphicx}
\usepackage{rotating}
\usepackage{fixmath}
\usepackage[pdftex]{hyperref}
\usepackage{soul}
\usepackage{stmaryrd}
\usepackage{relsize}
\usepackage{amssymb}

\hypersetup{
 colorlinks = true,
 urlcolor = blue,
 linkcolor = red,
}

\newcommand{\pe}[2]{\left\langle#1,#2\right\rangle}
\newcommand{\racion}[2]{\mbox{\small$\frac{{#1}}{{#2}}$}}

\newtheorem{thm}{Theorem}[section]
\newtheorem{cor}[thm]{Corollary} 
\newtheorem{lem}[thm]{Lemma}

\newtheorem{defn}[thm]{Definition}
\newtheorem{rem}[thm]{Remark}

\numberwithin{equation}{section}

\definecolor{Mycolor2}{HTML}{e85d04}
\sethlcolor{black}
 
\newcommand\poro[1]{{\textcolor{black}{#1}}} 
\begin{document}
\renewcommand{\PaperNumber}{***}

\FirstPageHeading

\ShortArticleName{Integral representations for a generalized Hermite linear functional}

%%%%%%%%%%%%%%%%%%%%%%%%%%%%%%%%%%%%%%%%%%%%%%%%%%%%
\ArticleName{Integral representations for a generalized Hermite linear functional}

\Author{
Roberto S. Costas-Santos$\,^{\ast}$
}
\AuthorNameForHeading{R.~S.~Costas-Santos
}
\Address{$^\ast$ Department of Quantitative Methods, Universidad Loyola Andaluc\'ia, E-41704, Sevilla, Spain
} 
% Address of First Author
\URLaddressD{
\href{http://www.rscosan.com}
{http://www.rscosan.com}
}
\EmailD{rscosa@gmail.com} 
\ArticleDates{Received~\today~in final form ????; Published online ????}
%%%%%%%%%%%%%%%%%%%%%%%%%%%%%%%%%%%%%%%%%%%%%%%%%%%%

% Abstract (Do not insert blank lines, i.e. \\) 
\Abstract{In this paper, we find new integral representations
for the \poro{generalized Hermite linear functional}
in the real line and the complex plane.
As an application, new integral representations for the Euler Gamma function are given.}

% Keywords

% The fields PACS, MSC, and JEL may be left empty or commented out if not applicable
%\PACS{J0101}
\Classification{33C45; 42C05 (Primary); 30E20; 33B15 (Secondary)}
%\JEL{}

%%%%%%%%%%%%%%%%%%%%%%%%%%%%%%%%%%%%%%%%%%
% Only for the journal Diversity
%\LSID{\url{http://}}

%%%%%%%%%%%%%%%%%%%%%%%%%%%%%%%%%%%%%%%%%%
% Only for the journal Applied Sciences:
%\featuredapplication{Authors are encouraged to provide a concise description of the specific application or a potential application of the work. This section is not mandatory.}
%%%%%%%%%%%%%%%%%%%%%%%%%%%%%%%%%%%%%%%%%%

%%%%%%%%%%%%%%%%%%%%%%%%%%%%%%%%%%%%%%%%%%
% Only for the journal Data:
%\dataset{DOI number or link to the deposited data set in cases where the data set is published or set to be published separately. If the data set is submitted and will be published as a supplement to this paper in the journal Data, this field will be filled by the editors of the journal. In this case, please make sure to submit the data set as a supplement when entering your manuscript into our manuscript editorial system.}

%\datasetlicense{license under which the data set is made available (CC0, CC-BY, CC-BY-SA, CC-BY-NC, etc.)}

%%%%%%%%%%%%%%%%%%%%%%%%%%%%%%%%%%%%%%%%%%
% Only for the journal Toxins
%\keycontribution{The breakthroughs or highlights of the manuscript. Authors can write one or two sentences to describe the most important part of the paper.}

%\setcounter{secnumdepth}{4}
%%%%%%%%%%%%%%%%%%%%%%%%%%%%%%%%%%%%%%%%%%
\section{Introduction} %\label{sec1}
The integral representation of special functions provides an alternative way 
to express these functions in terms of integrals involving other functions. 
They often involve a weight function and a kernel function related to the specific special 
function being considered. The weight function appears as a factor in the integral and 
reflects the orthogonality property of the associated orthogonal polynomials, and the 
kernel function represents the additional dependence.

The integral representation allows us to express special functions as infinite series or 
integrals involving some classical orthogonal polynomials. This connection arises 
from the fact that the orthogonality condition is satisfied by classical orthogonal polynomials, which
naturally leads to the appearance of these polynomials in the integral representation 
of special functions. In this work, we are going to consider the Hermite polynomials.

The hypergeometric functions, which have applications in many areas, including 
mathematical physics and combinatorics, can be represented in terms of integrals 
involving other hypergeometric functions and classical orthogonal polynomials like 
the Jacobi, Hermite, and Laguerre polynomials, which can be 
expressed as hypergeometric series (see \mbox{c.f. \cite{andrewsetal} and \cite{dlmf}({Section 16}})).

For a detailed history of the subject of integral representations 
for hypergeometric series and basic hypergeometric functions (which is a natural extension of 
the hypergeometric series),  see \cite{gas1989} and \cite{gasrah} ({Chapter 4}).

R. Sfaxi has established in \cite{sfa1}, by
means of a linear isomorphism, the so-called
{\it intertwining operator} on polynomials, a relationship between the ordinary 
Hermite polynomials and their analog nonsingular and of Laguerre--Hahn with
class zero.
Among others, the author has put in value an important linear functional, namely {\it the 
generalized Hermite linear functional}, denoted by $\mathscr{G}_{H}(\tau)$
of index $\tau \in \mathbb{C}$, with $\tau\ne -n,\; n\geq 1$, where their moments
are given by
\begin{equation} \label{1}
{\big(}\mathscr{G}_{H}(\tau){\big)}_{n}:=
\pe {\mathscr{G}_{H}(\tau)}{x^{n}}=
\left\{
\begin{array}{r@{\ \rm if \ }l} \displaystyle \dfrac{(\tau+1)_{2k}}
{k!2^{2k}}, & n=2k,\\[3mm] 0,& n=2k+1, \end{array}
\right.
\end{equation}
where $(a)_n$ is the Pochhammer symbol, defined as
\[
(a)_0:=1,\quad (a)_k:=a(a+1)\cdots (a+k-1),\quad
a\in \mathbb{C}\setminus\{0\}, \ k=1,2,3,\dots,
\]
thus $\mathscr{G}_{H}(\tau)$ is symmetric and monic,
i.e., ${\big(}\mathscr{G}_{H}(\tau){\big)}_{0}=1$.

Observe that setting $\tau=0$ in (\ref{1}) we recover
the Hermite linear functional, i.e., \mbox{$\mathscr{G}_{H}
\equiv \mathscr{G}_{H}(0)$}, that is well-known by
its integral representation
\begin{equation}\label{1.2}
\langle \mathscr{G}_{H},p\rangle =\frac{1}
{\sqrt{\pi}}\int_{-\infty}^{\infty}p(x)e^{-x^2}
dx,\quad p\in \mathbb{P}.
\end{equation}

\poro{So we can write}
\[ %\begin{equation}\label{1.2.1}
{\big(}\mathscr{G}_{H}(\tau){\big)}_{n}=\dfrac{(
\tau+1)_n}{(1)_n}(\mathscr{G}_{H})_n,\quad n=0, 1, \dots
\] %\end{equation}

\poro{Note that }the linear functional $\mathscr{G}_{H}$
is classical, since it is quasi-definite and satisfies
the Pearson equation
\begin{equation}\label{1.3}
 \mathscr{G}_{H}'+2x\mathscr{G}_{H}=0.
\end{equation}

\poro{Taking this} into account, the following result holds.

\begin{lem}
For any $\tau \in \mathbb{C}$, the linear functional
$\mathscr{G}_{H}(\tau)$ fulfills the difference equation
\[
{\big(}x^2\mathscr{G}_{H}(\tau){\big)}''+{\big(}2x
(x^2-\tau-2)\mathscr{G}_{H}(\tau){\big)}'+ {\big(}-
4x^2+(\tau+1)(\tau+2){\big)}\mathscr{G}_{H}(\tau)=
0.
\]
\end{lem}
\begin{proof}
Let $\tau \in \mathbb{C}$, if we define the linear functional
$\mathscr{E}(\tau)$ as
\[
\mathscr{E}(\tau):={\big(}x^2\mathscr{G}_{H}(\tau)
{\big)}''+{\big(}2x(x^2-\tau-2)\mathscr{G}_{H}(\tau)
{\big)}'+{\big(}-4x^2+(\tau+1)(\tau+2){\big)}
\mathscr{G}_{H}(\tau).
\]

\poro{Then}, for $n\ge 0$, one obtains
\begin{equation}\label{1.5}
{\big(}\mathscr{E}(\tau){\big)}_{n}=-2(n+2){\big(}
\mathscr{G}_{H}(\tau){\big)}_{n+2}+(n+\tau+2)(n+
\tau+1){\big(}\mathscr{G}_{H}(\tau){\big)}_{n}.
\end{equation}

\poro{Since} $\mathscr{G}_{H}(\tau)$ is
symmetric, then ${\big(}\mathscr{E}(\tau){\big)
}_{2k+1}=0$, for every $k\geq 0$.
On the other hand, setting $n=2k$ in (\ref{1.5}) and
taking into account
(\ref{1}), we get for $k\ge 0$,
\[\begin{array}{rl@{\!\!}}
{\big(}\mathscr{E}(\tau){\big)}_{2k}&=\;-4(k+1){\big(}
\mathscr{G}_{H}(\tau){\big)}_{2k+2}+(2k+\tau+2)
(2k+\tau+1){\big(}\mathscr{G}_{H}(\tau){\big)}_{2k}\\[2mm]
&=\;- \dfrac{(\tau+1)_{2k+2}}{k!2^{2k}}+
\dfrac{(2k+1+\tau+1)(2k+\tau+1)
(\tau+1)_{2k}}{k!2^{2k}}\\[2mm]& = 0.
\end{array}
\]

\poro{Therefore,} ${\big(}\mathscr{E}(\tau)
{\big)}_{n}=0$ for all $n=0, 1, \dots$. Hence, the result holds.
\end{proof}

Our purpose in this work is to provide integral
representations for the linear functional $\mathscr{
G}_{H}(\tau)$, either on the real axis, or on the
complex plane. %AUTHOR: I removed the newline command
\poro{More} precisely, the problem consists of determining
a weight function $G_{H}(\bullet;\tau)$, such that
\[
\langle \mathscr{G}_{H}(\tau),p\rangle
=\int_{\Omega}p(x)G_{H}(x;\tau)dx,\quad p\in
\mathbb{P},
\]
where $\Omega$ is an interval in the real line, or 
a contour in the complex plane.

The paper is organized as follows.
In the next section, there are some preliminaries and notations.
In Sections \ref{sec3} and \ref{sec4}, integral representations in the
real line and in the complex plane, respectively, are
provided.
As an application of the previous results, in Section \ref{sec5}, some 
new integral representations for the Euler Gamma function are given.

\section{Preliminaries and Notation} \label{sec2}
Let $\mathbb{P}$ be the vector space of polynomials
with complex coefficients and let $\mathbb{P}'$ be
its dual space.
We denote by $\langle u,f\rangle$ the action of the
linear functional $u\in {\mathbb{P}'}$ on the polynomial
$f\in{\mathbb{P}}$.
In particular, we denote by $(u)_{n}:=\langle u,x^n
\rangle$, $n\geq 0$, the moments of $u$.
%%Moreover, for any linear functional $u$, and any polynomial
%%$f$, let $u'$, $fu$ be the linear functionals
%%defined by duality, as
%%$$\begin{array}{rl}
%%\langle u',p\rangle&:=\;-\langle u, p'\rangle,\quad
%% p\in \mathbb{P},\\
%%\langle fu,p\rangle&:=\;\langle u, fp\rangle, \quad
%% p\in \mathbb{P}.
%%\end{array}$$

\begin{defn} A linear functional $u$ is called 
{\poro{symmetric}} 
if $(u)_{2n+1}=0$, for 
all $n=0, 1, \dots$, and it is called 
\poro{monic} if $(u)_0=1$.
\end{defn}

{In fact, for any $\tau \in \mathbb{C}$, the linear functional 
$\mathscr{G}_{H} (\tau)$ is symmetric (see \eqref{1})
which allows us to suppose the weight function 
$G_{H} (\bullet;\tau)$ is even, i.e., it can be written as 
\mbox{$G_{H}(x;\tau)=U(|x|;\tau)$}}, where $U(\bullet;\tau)$ is a
function defined on $(0,\infty)$.
In fact, this is a direct consequence of the following
result.
\begin{lem} %\label{lem2.1}
Let $\mathscr{L}$ be a symmetric linear function
that has an integral representation.
Then, there exists a function $U$ defined on
$(0,\infty)$, such that
\[
\pe {\mathscr{L}}p =\int_{-\infty}^{\infty}p(x)
U(|x|)dx.
\]
\end{lem}
\begin{proof}
From the assumption there exists a function $L$, defined
on $(-\infty,\infty)$, such that
\[
\pe {\mathscr{L}}p =\int_{-\infty}^{\infty}
p(x) L(x)dx.
\]

\poro{Let us} introduce the following two functions, defined on
$(0,\infty)$, as follows:
\[
U(x)=\dfrac{L(x)+L(-x)}{2},\quad V(x)=\left\{
\begin{array}{rl}
\dfrac{L(x)-L(-x)}{2x}, & \;{\rm if}\;x\ne 0, \\
0, & \;{\rm if}\;x=0. \end{array} \right.
\]

\poro{A straightforward} calculation gives that $L(x)=U(|x|)+
xV(|x|)$, for all $x\in \mathbb R$. Moreover, 
since $x^{2n+1} V(|x|)$ is an odd function we have
\[%\begin{array}{rl}
(\mathscr{L})_{2n}=\;\int_{-\infty}^{\infty}x^{2n}
U(|x|)dx+\int_{-\infty}^{\infty}x^{2n+1} V(|x| )dx=
\int_{-\infty}^{\infty}x^{2n} U(|x|)dx.
%\end{array}
\]

\poro{On the other} hand, since $\mathscr{L}$ is symmetric
and $x^{2n+1} U(|x| )$ is an odd function, we get
\[
(\mathscr{L})_{2n+1}=\;\int_{-\infty}^{\infty}
x^{2n+1} U(|x|)dx=0.
\]

\poro{Therefore}, for any polynomial $p\in \mathbb{P}$,
\[
\pe {\mathscr{L}} p=\;\int_{-\infty}^{\infty}p(x)
 U(|x|)dx.
\]
\end{proof}

The next result related to hypergeometric
functions will be useful later.
\begin{lem}[\cite{kolest,leb}]  %{\cite{kolest,leb}}  %\label{lemma-2.2}
The following formulae hold:
\begin{enumerate}
\item If $\Re(\alpha)>0$ and $\Re(s)>0$,
then
\begin{equation}\label{2.1}
\int_{0}^{\infty}t^{\alpha-1}{_1F_1}(a_1;b_1;t)
e^{-st}dt=\dfrac{\Gamma(\alpha)}{s^{\alpha}}\,
{_2F_1}(a_1,\alpha;b_1;1/s).
\end{equation}
 \item If $\Re(c-a-b)>0$, then
\begin{equation}\label{2.2}
{_2F_1}(a,b;c;1)=\dfrac{\Gamma(c)\Gamma(c-a-b)}
{\Gamma(c-a)\Gamma(c-b)},
\end{equation}
\end{enumerate}
where 
\[
{}_2F_1(a,b;c;z):=\sum_{k=0}^\infty \dfrac{(a)_k(b)_k}{(c)_k}\dfrac{z^k}{k!},\qquad 
{}_1F_1(a;b;z):=\sum_{k=0}^\infty \dfrac{(a)_k}{(b)_k}\dfrac{z^k}{k!}.
\]
\end{lem}
In future work, we will denote by $H_\tau(x)$ the {\it \poro{Hermite
function (of degree} $\tau$)},
which can be represented in terms of the confluent
hypergeometric function ${{}_1F_1}$ as follows %Please check that intended meaning has been retained.  
%AUTHOR: Checcked
\cite{leb}:
\begin{equation}\label{2.3}
H_\tau(x)=2^\tau\frac{\Gamma(\frac{1}{2})}{\Gamma(
\frac{1-\tau}{2})}\,{{}_1F_1}\left(-\frac{\tau}{2};\frac 12;
x^2\right)+2^{\tau}x\frac{\Gamma(-\frac{1}{2})}
{\Gamma(-\frac{\tau}{2})}\,{{}_1F_1}\left(\frac{1-\tau}2;
\frac{3}{2};x^2\right).
\end{equation}

\section{Integral Representation on $\mathbb{R}$} \label{sec3}
In the following result, we present a new definite
integration formulae involving the Hermite functions.
\begin{lem} %\label{prop-3.1}
For any $(z,\tau)\in \mathbb{C}^2$, with $\Re(z)>-1$,
the following formulae hold:
\begin{equation}\label{3.1}
\int_{0}^{\infty}x^{z}H_\tau(x)e^{-x^2}dx=\dfrac{
\sqrt{\pi}}{2^{z-\tau+1}}\dfrac{\Gamma(z+1)}{\Gamma(
\frac{z-\tau}{2}+1)},
\end{equation}
\begin{equation}\label{3.2}
\int_{-\infty}^{\infty}|x|^{z}H_\tau(|x|)e^{-x^2}dx=
\dfrac{\sqrt{\pi}}{2^{z-\tau}}\dfrac{\Gamma(z+1)}
{\Gamma(\frac{z-\tau}{2}+1)}.
\end{equation}
%where $H_\tau(x)$ is the Hermite function (of degree $\tau$).
\end{lem}

\begin{proof}
Since the function $|x|^\nu H_\tau(|x|) e^{-x^2}$
is even, it is enough to prove (\ref{3.1}).

Let us fix $\tau\in \mathbb{C}$, with $\Re(\tau)>-1$.
For any $z\in \mathbb{C}$, with $-1<\Re(z)<\Re(\tau)$,
let us consider the following integral:
\[
\Lambda(z):=\int_{0}^{\infty}x^{z}H_\tau(x)e^{-x^2}dx.
\]

\poro{Using} (\ref{2.3}), the previous integral can be written
as
\begin{equation}\label{3.3}
\Lambda(z)=2^\tau\frac{\Gamma(\frac{1}{2})}
{\Gamma(\frac{1-\tau}{2})}\Pi(z)+ 2^\tau
\frac{\Gamma(-\frac{1}{2})}{\Gamma(-\frac{\tau}{2})}
\Omega(z),
\end{equation}
where
\[\begin{array}{rl}
\Pi(z)&:=\displaystyle \int_{0}^{\infty}x^{z}{_1F_1}
\left(-\frac{\tau}{2};\frac{1}{2};x^2\right)e^{-x^2}dx,
\\[5mm] \Omega(z)&:=\displaystyle \int_{0}^{\infty}x^{z+1}{_1F_1}
\left(\frac{1-\tau}{2}; \frac{3}{2};x^2\right)e^{-x^2}dx.
\end{array}\]

\poro{By changing} the variable of integration, by setting
$t=x^2$, and using (\ref{2.1}), with $s=1,$ $\alpha=
(z+1)/2$, $a_1=-\tau/2$, and $b_1=1/2$,
we obtain
\[
\Pi(z)=\frac{1}{2}\,\Gamma\left(\frac{z+1}{2}\right){_2F_1}
\left(-\frac{\tau}{2},\frac{z+1}{2};\frac 12;1\right).
\]

\poro{Again}, with (\ref{2.1}), where $s=1$, $\alpha=
(z+2)/2$, $a_1=(1-\tau)/2$, and $b_1=3/2$,
we get
\[
\Omega(z)=\frac{1}{2}\, \Gamma\left(\frac{z+2}{2}\right)
{_2F_1}\left(\frac{1-\tau}{2},\frac{z+2}{2};\frac 32;1\right).
\]

\poro{Since} $\Re(z)<\Re(\tau)$, by using (\ref{2.2})
$\Pi(z)$ and $\Omega(z)$ can be written as
\[
\begin{array}{rl}
\Pi(z)&= \dfrac{\Gamma(\frac{z+1}{2})
\Gamma(\frac{1}{2})\Gamma(\frac{\tau-z}{2})}
{2\Gamma(\frac{1+\tau}{2}) \Gamma(-\frac{z}{2})},\\[4mm]
\Omega(z)&=\dfrac{\Gamma(\frac{z+2}{2})
\Gamma(\frac{3}{2})\Gamma(\frac{\tau-z}{2})}{2\Gamma(
\frac{2+\tau}{2})\Gamma(\frac{1-z}{2})}.
\end{array}
\]

\poro{Therefore}, taking into account 
%\[
$\Gamma(\racion{1}{2})^2=-\Gamma(-\racion{1}{2})
\Gamma(\racion 32)=\pi$,
%\]
the expression (\ref{3.3}) can be rewritten as follows:
\[
%\begin{equation}\label{3.4}
\Lambda(z)=\dfrac{2^{\tau-1}\pi\Gamma(
\frac{\tau-z}{2})}{\Gamma(-\frac{z}{2})\Gamma(\frac{1-
z}{2})}{\Big(}U(z,\tau)-U(z+1,\tau+1){\Big)},
\]%\end{equation}
where
\[
U(z,\tau)=\dfrac{\Gamma(\frac{z+1}{2})\Gamma(\frac{1-z}2)}
{\Gamma(\frac{1+\tau}{2})\Gamma(\frac{1-\tau}{2})}.
\]

\poro{Using} the duplication formula
\[
\Gamma(u)\Gamma(1-u)=\dfrac{\pi}{\sin(\pi
u)},
\]
a straightforward calculation leads to
\[
%\begin{equation}\label{3.19}
U(z,\tau)=\frac{\cos(\frac{\pi}{2}\tau)}{\cos(\frac{\pi}
{2}z)},\quad U(z+1,\tau+1)=\frac{\sin(\frac{\pi}{2}\tau)}
{\sin(\frac{\pi}{2}z)}.
\]%\end{equation}

\poro{Then},
\[
\Lambda(z)=-\frac{2^{\tau}\pi\Gamma(\frac{\tau-z}{2})}
{\Gamma(-\frac{z}{2})\Gamma(\frac{1-z}{2}) }
\frac{\sin{\big(}\frac{\pi}{2}(\tau-z){\big)}}{\sin(\pi
z)},
\]
so, by using the Gauss--Legendre multiplication formula,
\[
\Gamma(u)\Gamma(u+\racion{1}{2})= 2^{1-2u}\sqrt{\pi}\,
\Gamma(2u),
\]
and, again, with the duplication formula, we get
\[
\Lambda(z)=\frac{\sqrt{\pi}}{2^{z-\tau+1}}
\frac{\Gamma(z+1)}{\Gamma(1+\frac{z-\tau}{2})}.
\]

\poro{For this} proof, we assumed the conditions
$-1<\Re(z)<\Re(\tau)$, then the integral $\Lambda(z)$ 
converged exponentially to zero when $\tau\to \infty$.
Hence, through analytic continuation, (\ref{3.3})
is valid for each $(\tau,z)\in \mathbb{C}^2,$
with $\Re(z)>-1.$
\end{proof}
\begin{rem}
Note that the above result also covers the $z=\tau$ case.
In fact, if $\tau=$0, 1, $\dots$ this identity
represents the property of orthogonality for the monic
Hermite polynomials.
\end{rem}
As a consequence, we have the following result:
\begin{cor} \label{cor3.3}
For any $\tau\in \mathbb{C}$, with $\Re(\tau)>-1$,
the following formulae hold:
\begin{equation}\label{3.4}
\int_{0}^{\infty}x^{2n+\tau}H_\tau(x)e^{-x^2}dx=
\frac{\sqrt{\pi}}{2^{2n+1}}\frac{\Gamma(2n+\tau+1)}
{\Gamma(n+1)},
\end{equation}
\begin{equation}\label{3.5}
\int_{-\infty}^{\infty}x^{2n}|x|^{\tau}H_\tau(|x|)
e^{-x^2}dx=\frac{\sqrt{\pi}}{2^{2n}}\frac{\Gamma(2n+
\tau+1)}{\Gamma(n+1)}.
\end{equation}
\end{cor}
\begin{thm}
For any $\tau\in \mathbb{C}$, with $\Re(\tau)>-1$, 
the linear functional $\mathscr{G}_{H}(\tau)$ has
the following integral representation:
\begin{equation} \label{3.16}
\langle
\mathscr{G}_{H}(\tau),p\rangle =\frac{1}{\sqrt{\pi}
\, \Gamma(\tau+1)}\int_{-\infty}^{\infty} p(x)|x|^\tau
H_\tau(|x|) e^{-x^2}dx,\quad p\in\mathbb{P},
\end{equation}
where $H_\tau$ is the Hermite function (of degree $\tau$).
\end{thm}
\begin{proof}
Due to the Equation \eqref{1} and Corollary \ref{cor3.3}, 
\[\begin{array}{rl}
{\big(}\mathscr{G}_{H}(\tau){\big)}_{2n}&
=\displaystyle \frac{(\tau+1)_{2n}}{n!2^{2n}}=\frac{\Gamma(2n+
\tau+1)}{2^{2n}\Gamma(n+1)\Gamma(\tau+1)}\\[4mm]
&=\displaystyle \frac{1}{\sqrt{\pi}\,\Gamma(\tau+1)}
\int_{-\infty}^{\infty}x^{2n}|x|^{\tau}H_\tau(|x|)
e^{-x^2}dx,\\[4mm] {\big(}\mathscr{G}_{H}(\tau){\big)}_{2n+1}&
=0=\displaystyle \frac{1}{\sqrt{\pi}\,\Gamma(\tau+1)}
\int_{-\infty}^{\infty} x^{2n+1}|x|^{\tau}H_\tau(|x|)
e^{-x^2}dx.
\end{array}
\]
\poro{Therefore}, one has
\[
{\big(}\mathscr{G}_{H}(\tau){\big)}_{n}=\frac{1}
{\sqrt{\pi}\,\Gamma(\tau+1)}\int_{-\infty}^{\infty}x^{n}
|x|^{\tau}H_\tau(|x|)e^{-x^2}dx,\qquad n=0, 1, \dots
\]
\poro{Consequently}, for any polynomial $p\in
\mathbb{P}$,
\[
\pe {\mathscr{G}_{H}(\tau)}{p}=
\frac{1}{\sqrt{\pi}\,\Gamma(\tau+1)}
\int_{-\infty}^{\infty}p(x)|x|^{\tau}
H_\tau(|x|)e^{-x^2}dx.
\]
\end{proof}

Observe that if we set $n=0$ in \eqref{3.4}
we get a new integral representation for the
Euler Gamma function. In fact, for
any $\tau\in \mathbb{C}$, with $\Re(\tau)>-1$,
\begin{equation}\label{3.6}
\Gamma(\tau+1)=\frac{2}{\sqrt{\pi}}\int_{0}^{\infty}
x^{\tau}H_\tau(x)e^{-x^2}dx,
\end{equation}
\begin{equation}\label{3.7}
\Gamma(\tau+1)=\frac{1}{\sqrt{\pi}}\int_{-\infty}^{\infty}
|x|^{\tau}H_\tau(|x|)e^{-x^2}dx.
\end{equation}

\section{Integral Representation on the Complex Plane} \label{sec4}
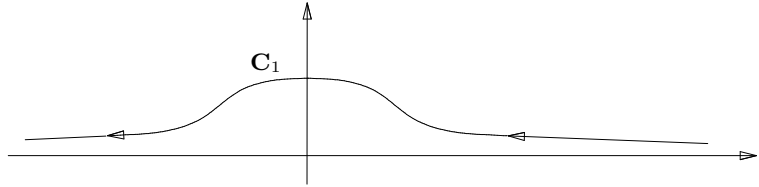
\begin{figure}[!hbt]
\begin{tikzpicture}[domain=-7:7, scale=0.62]
\draw[gray]
plot[domain=-7:-4] (\x, -0.67/\x);
\draw[gray,latex-] 
plot[domain=-4:-1.2] (\x, -0.67/\x) ;
\draw[gray] 
plot[domain=0:1.2] (\x,{2*(1-\x^2+ 0.5*\x^4-0.167*\x^6+0.042*\x^8)}) 
plot[domain=-1.2:0] (\x,{2*(1+\x^2-0.5*\x^4+0.167*\x^6-0.042*\x^8)})
node[right,black]{\tiny $\bf C_1$};
\draw[gray,latex-] 
plot[domain=1.2:7] (\x, 0.67/\x);
\draw[->] (-7,0) -- (7.3,0) node[right] {$x$};
\draw[->] (0,-0.25) -- (0,3) node[left] {$y$};
\end{tikzpicture}
%\end{center}
\caption{\label{gr1}\poro{Path ${\bf C_1}$}}
\end{figure}

\begin{thm} For any $\tau\in \mathbb{C}$, the following identities hold:
\begin{enumerate}
\item[(i)] 
\[
\int_{\bf{C}_1}\zeta^{2n+1}|\zeta|^{\tau}
H_\tau(|\zeta|)e^{-\zeta^2}d\zeta=0,\quad n=0, 1, \dots
\]
\item[(ii)] For any $n\in \mathbb N$, so that $\tau+2n+1$ is not a negative 
integer, we have 
\[
\int_{\bf{C}_1}\zeta^{2n}|\zeta|^{\tau}
H_\tau(|\zeta|)e^{-\zeta^2}d\zeta=-\dfrac{\sqrt{\pi}}{2^{2n}} 
\dfrac{\Gamma(2n+\tau+1)}{\Gamma(n+1)},\quad n=0, 1, \dots
\]
\end{enumerate} 
where $\bf{C}_1$ is the following contour in the
\poro{complex plane (See Figure} \ref{gr1}).%mdpi: Please check if this should be formatted as Figure 1. xxx. If so, please add the caption and its citation in the text accordingly. Same as this issue below.
%AUTHOR: Amended
 \end{thm}

%%%%%%%%%%%%%%%%%%%%%% Begin figure 1 %%%%%%%%
%\begin{center}
%%%%%%%%%%%%%%%%%%%%%% End figure 1 %%%%%%%%%
\begin{proof} We deform $\bf{C}_1$ into a contour
$\tilde{\bf{C}}_1$ consisting of
two straight lines and a circle \poro{(see Figure} \ref{gr2}).

%%%%%%%%%%%%%%%%%%%%%% Begin figure 2 %%%%%%%%%%
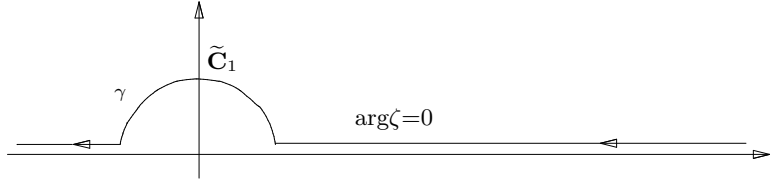
\begin{figure}[!ht]
%\begin{center}
\begin{tikzpicture}[domain=-7:7, scale=0.62, samples=200]
\draw[gray] 
plot[domain=-7:-3] (\x, 0.2);
\draw[gray,latex-] 
plot[domain=-3:-2.1] (\x,0.2) -- plot[domain=pi-0.4:0.4] ({2*cos(\x r)}, {2*sin(\x r)}) 
-- plot[domain=2.1:4] (\x,0.2);   
\draw  (0.25, 2.2) node[right,black]{\tiny $\bf \widetilde C_1$}
(4,0.2) node[above,black] {\tiny arg($\zeta$)=0}
(-1.8,1) node[above]{\tiny $\gamma$}; 
\draw[gray,latex-] 
plot[domain=4:7] (\x,0.2);
\draw[->] (-7,0) -- (7.3,0) node[right] {$x$};
\draw[->] (0,-0.25) -- (0,3) node[left] {$y$};
\end{tikzpicture}
\caption{\label{gr2}\poro{Path $\bf \widetilde C_1$}}
%\end{center}
\end{figure}
%%%%%%%%%%%%%%%%%%%%%% End figure 2 %%%%%%%%%%%

\noindent
where $\gamma:=\{\zeta \in \mathbb{C}: \Im(\zeta)>0,
|\zeta|=\epsilon\}$, being $\epsilon>0$.

Now, for each integer $n\geq 0$ and $\tau\in \mathbb{C},$
we define
\[
\begin{array}{rl}
I_n(\tau):=& \displaystyle \int_{\tilde{\bf{C}}_1}\zeta^n| \zeta
|^{\tau}H_\tau(|\zeta|)e^{-\zeta^2}d\zeta=\int_{\infty}^{\epsilon}\!\!\zeta^n| \zeta
|^{\tau}H_\tau(|\zeta|)e^{-\zeta^2}d\zeta \\[3mm]
&\displaystyle+\int_{\gamma}
\!\!\zeta^n|\zeta|^{\tau}H_\tau(|\zeta|)e^{-\zeta^2}d\zeta
+\int_{-\epsilon}^{-\infty} \!\!\!\!
\zeta^n|\zeta|^{\tau}H_\tau(|\zeta|)e^{-\zeta^2}d\zeta.
\end{array} 
\]

\poro{So}, if $\Re(\tau)>-n-1$, after a direct computation, we get
\[
\begin{array}{rl}
\displaystyle \lim_{\epsilon\rightarrow
0}\int_{\infty}^{\epsilon}\zeta^n|\zeta
|^{\tau}H_\tau(|\zeta|)e^{-\zeta^2}d\zeta&=
\displaystyle -\int_{0}^{\infty}x^{n+\tau}
H_\tau(x)e^{-x^2}dx,\\ \displaystyle \lim_{
\epsilon\rightarrow 0}\int_{-\epsilon}^{-\infty}
\zeta^n|\zeta|^{\tau}H_\tau(|\zeta|)e^{-\zeta^2}
d\zeta&=\displaystyle -(-1)^n\int_{0}^{\infty}
x^{n+\tau}H_\tau(x)e^{-x^2}dx.
\end{array}
\]

\poro{For} the middle integral, we obtain
\[
\begin{array}{rl}
\displaystyle {\Big|} \int_{\gamma}\zeta^n|\zeta|^{\tau}
H_\tau(|z|)e^{-\zeta^2}d\zeta {\Big|}
&\displaystyle =\;{\Big|} \int_{0}^{\pi}\epsilon^n e^{in\theta}
\epsilon^\tau H_\tau(\epsilon )e^{-\epsilon^2 e^{2i
\theta}}\epsilon i e^{i\theta}d\theta{\Big|}\\
&\displaystyle \leq\epsilon^{n+\Re(\tau)+1} \int_{0}^{\pi}|H_\tau(
\epsilon )| e^{-\epsilon^2 \cos(2\theta)}d\theta,
\end{array}
\]
knowing that $H_\tau(0)=2^\tau\sqrt{\pi}/\Gamma(\frac{1-\tau}{2})$, it
is straightforward to see that
\[
\lim_{\epsilon\to 0} \int_{\gamma}\zeta^n|
\zeta|^{\tau}H_\tau(|\zeta|)e^{-\zeta^2}d\zeta =\;0.
\]

\poro{Therefore}, for each $n\geq 0$ and $\tau\in
\mathbb{C}$, such that $\Re(\tau)>-n-1$, we have
\[
I_n(\tau)=-{\big(}(-1)^n+1{\big)}\int_{0}^{\infty}x^{n+
\tau }H_\tau(x)e^{-x^2}dx.
\]

\poro{Then,} $I_{2n+1}(\tau)=0$ for all $n\ge 0$.
Notice that for the proof of i), we assumed $\Re(\tau)>-n-1$,
but the integral converges exponentially when $\tau \to
\infty$, and therefore it exists for all $\tau$.
Hence, (i) holds through analytic continuation for
any $\tau\in \mathbb{C}$.

{On the other }hand, using (\ref{3.4}), it follows that
\[
I_{2n}(\tau)=-\frac{\sqrt{\pi}}{2^{2n}}
\frac{\Gamma(2n+\tau+1)}{\Gamma(n+1)}.
\]

\poro{Hence, (ii) holds,} for the same reason already quoted and
by analytic continuation of $\tau \in \mathbb{C}$,
except when $2n+\tau+1$ is a negative integer, where the
function $\Gamma$ is undefined.
\end{proof}
As a consequence, we have the following result.
\begin{thm}
For any $\tau\in \mathbb{C}$, with $-\tau\not\in \mathbb{N}$, 
the linear functional $\mathscr{G}_{H}(\tau)$ has the following 
integral representation:
\begin{equation} \label{3.16p}
\langle \mathscr{G}_{H}(\tau),p\rangle=-\frac{1}
{\sqrt{\pi}\Gamma(\tau+1)}\int_{\bf{C}_1} p(x)|x|^\tau
H_\tau(|x|)e^{-x^2}dx,\quad p\in\mathbb{P},
\end{equation}
where $H_\tau$ is the Hermite function (of degree $\tau$).
\end{thm}
Using an analog idea allows us to formulate another
integral representation for the gamma function in
the complex plane by using a different contour.
\begin{thm} %\label{thm4.3}
For any $\tau\in \mathbb{C}$, with $-\tau\not\in \mathbb{N}$, 
the Euler's  Gamma function satisfies the following integral 
representation:
\begin{equation} \label{3.21}
\Gamma(\tau+1)=\frac{2}{\sqrt{\pi}(e^{2\pi i\tau}-1)}
\int_{\bf{C}}\zeta^{\tau}H_\tau(\zeta)e^{-\zeta^2}d\zeta,
\end{equation}
where \poro{$C$}%mdpi: Please confirm if the bold is unnecessary and can be removed
%AUTHOR: style \bf removed 
is the following contour in the \poro{complex
plane} 

%%%%%%%%%%%%%%%%%%%%%% Begin figure 3 %%%%%%%%%%
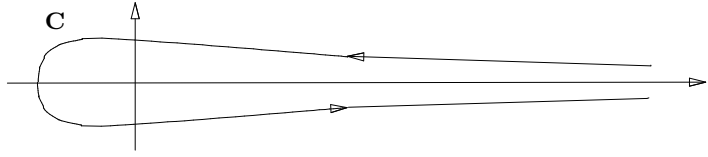
\begin{figure}[!hbt]
%\begin{center}
\begin{tikzpicture}[domain=-2:6, scale=0.62]
\draw[gray]
plot[domain=6:4] (\x, 0.67/\x)
plot[domain=6:4] (\x, -0.67/\x);
\draw[gray, latex-] 
plot[domain=4:1.2] (\x, -0.67/\x);
\draw[gray, -latex] 
plot[domain=4:1.2] (\x, 0.67/\x);
\draw[gray] 
plot[domain=1.25:0] (\x,{2*(1-\x^2+0.5*\x^4-0.167*\x^6+0.042*\x^8)})
(-1.25,2.2) node[right,black]{\tiny $C$}
plot[domain=1.25:0] (\x,{-2*(1-\x^2+0.5*\x^4-0.167*\x^6+0.042*\x^8)});
\draw[gray] 
plot[domain=pi/2:3*pi/2] ({2*cos(\x r)}, {2*sin(\x r)});
%plot[domain=1.2:7] (\x, 0.67/\x);
\draw[->] (-3,0) -- (6.2,0) node[right] {$x$};
\draw[->] (0,-3) -- (0,3) node[left] {$y$};
\end{tikzpicture}
%\end{center}
\caption{\label{gr3}\poro{Path of integration $C$}}
\end{figure}
%%%%%%%%%%%%%%%%%%%%%% End figure 3 %%%%%%%%%%%
\end{thm}

\begin{proof} We deform \poro{${C}$} into a contour
\poro{$\tilde{{C}}$} consisting of two straight lines and
a circle{:} %\poro{as follows:}

%%%%%%%%%%%%%%%%%%%%%% Begin figure 4 %%%%%%%%%%
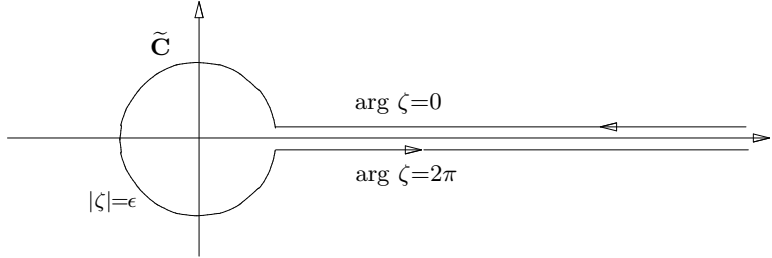
\begin{figure}[!hbt]
%\begin{center}
\begin{tikzpicture}[domain=-2:6, scale=0.62, samples=200]
\draw[gray] 
plot[domain=6:4] (\x, -0.2);
\draw[gray,latex-] 
plot[domain=4:2.1] (\x,-0.2) -- plot[domain=2*pi-0.4:0.4] ({2*cos(\x r)}, {2*sin(\x r)}) 
-- plot[domain=2.1:4] (\x,0.2);   
\draw  (-1.25, 2.2) node[right,black]{\tiny $\widetilde C$}
(4,0.2) node[above,black] {\tiny arg($\zeta$)=0}
(4,-0.2) node[below,black] {\tiny arg($\zeta$)=$2\pi$}
(-1.8,-2.2) node[above]{\tiny $|\zeta|=\epsilon$}; 
\draw[gray,latex-] 
plot[domain=4:6] (\x,0.2);
\draw[->] (-3,0) -- (6.2,0) node[right] {$x$};
\draw[->] (0,-3) -- (0,3) node[left] {$y$};
\end{tikzpicture}
%\end{center}
\caption{\label{g4}\poro{Path of integration $\widetilde C$}}
\end{figure}
%%%%%%%%%%%%%%%%%%%%%% End figure 4 %%%%%%%%%%%

\vskip 0.1cm

{We let}
\[
J(\tau)=\int_{\tilde{\bf{C}}}\zeta^{\tau}H_\tau(\zeta)
e^{-\zeta^2}d\zeta
\]
\noindent
Then
\[
J(\tau)=
\int_{\infty}^{\epsilon}\zeta^{\tau}H_\tau(\zeta)
e^{-\zeta^2}d\zeta+\int_{|\zeta|=\epsilon}\zeta^{\tau}
H_\tau(\zeta)e^{-\zeta^2}d\zeta +\int_{\epsilon}^{\infty}
\zeta^{\tau}H_\tau(\zeta)e^{-\zeta^2}d\zeta,
\]
and if $\Re(\tau)>-1$ in a direct way, we obtain 
\[
\begin{array}{rl}
\displaystyle \lim_{\epsilon\to 0}\int_{\infty}^{\epsilon}\zeta^{\tau}
H_\tau(\zeta)e^{-\zeta^2}d\zeta&=\displaystyle -\frac{\sqrt{\pi}}{2}
\Gamma(\tau+1),\\[6mm] \displaystyle \lim_{\epsilon\to 0}\int_{\epsilon}^{
\infty}\zeta^{\tau}H_\tau(\zeta)e^{-\zeta^2}d\zeta&\displaystyle
=e^{2\pi i\tau}\frac{\sqrt{\pi}}{2}\Gamma(\tau+1).
\end{array}
\]

\poro{For the} middle integral, we obtain
\[
\begin{array}{rl}
\displaystyle {\Big|} \int_{|\zeta|=\epsilon}\zeta^{\tau}
H_\tau(\zeta)e^{-\zeta^2}d\zeta|
&=\displaystyle \;\left| \int_{0}^{2\pi}(\epsilon
e^{i\theta})^{\tau}H_\tau(\epsilon e^{i\theta})
e^{-\epsilon^2 e^{2i\theta}}\epsilon i e^{i\theta}d\theta\right|\\
&\displaystyle \leq\epsilon^{\Re(\tau)+1} \int_{0}^{2\pi}| H_\tau(
\epsilon e^{i\theta})| e^{-\epsilon^2 \cos(2\theta)-
\theta{\big(}\Im(\tau)+1{\big)}}d\theta,
\end{array}
\]
thus,
\[
\lim_{\epsilon\to 0} \int_{|\zeta|=\epsilon}
\zeta^{\tau}H_\tau(\zeta)e^{-\zeta^2}d\zeta= 0.
\]

\poro{Finally,}
\[
J(\tau)=(e^{2\pi i\tau}-1)\frac{\sqrt{\pi}}{2}
\Gamma(\tau+1),
\]
hence, the result holds.
In the proof, we have assumed that $\Re(\tau)>-1$, but
the integral (\ref{3.21}) converges exponentially at
infinity, and therefore it exists for all $\tau$.
In fact, through analytic continuation, the result is valid
for every complex $\tau$, except for the negative
integers, where the denominator vanishes.
\end{proof}

In addition, from the last representation,
we obtain the following:
\[
\Gamma(\tau+1)=\frac{1}{i\sqrt{\pi}\sin(
\pi \tau)}\int_{\bf{C}}(-\zeta)^{\tau}H_\tau(\zeta)e^{-
\zeta^2}d\zeta.
\]

\poro{In the} last result, we show a representation for the reciprocal of
$\Gamma(\tau+1)$.
\begin{thm}
\[
\frac 1{\Gamma(\tau+1)}=-i\pi^{-\frac 32}\int_{\bf{C}}
(-\zeta)^{-1-\tau}H_{-1-\tau}(\zeta)e^{-\zeta^2}d\zeta.
\]

\poro{This representation} is valid for all $\tau$ and \poro{${C}$}
is the same contour as in the previous theorem.
\end{thm}
\begin{proof} Based on the last representation, one has
\[
\begin{array}{rl}
\Gamma(-\tau)=& \displaystyle \frac{1}{i\sqrt{\pi}\sin(\pi\tau)}
\int_{\bf{C}}(-\zeta)^{-1-\tau}H_{-1-\tau}(\zeta)
e^{-\zeta^2}d\zeta\\[5mm]=& \displaystyle \frac{\Gamma(\tau+1)
\Gamma(-\tau)}{i\pi^{\frac{3}{2}}}\int_{\bf{C}}
(-\zeta)^{-1-\tau}H_{-1-\tau}(\zeta)e^{-\zeta^2}d\zeta.
\end{array}
\]

\poro{This leads} to the desired result.
\end{proof}

\section{Conclusions}\label{sec5}
We have obtained  integral representations of a generalized linear Hermite functional,
which is among the natural extensions of the linear Hermite functional, using 
the fact this linear functional is symmetric, i.e., the odd moments associated with this 
functional are zero, and also the fact  that some hypergeometric representations associated 
with the Hermite polynomials are known. Observe that this can also be implemented for other 
symmetric classical orthogonal polynomials. Moreover, we have obtained an integral 
representation for the generalized linear Hermite functional in the complex plane, and 
from this integral representation, we are able to obtain a novel integral representation 
for the Euler Gamma function.

Of course, this method can be applied not only to other (symmetric) classical orthogonal 
polynomials but to any other symmetric orthogonal polynomial sequence for which 
a hypergeometric representation is known. This is something we should do 
in order to obtain novel integral representations for other Special functions; for example
we could consider some other generalization for the Hermite linear functional, 
as well as some Laguerre--Hahn or semi-classical, orthogonal polynomials 
(see, e.g., \cite{rebocho, cohletall2020} and the references therein).

%%%%%%%%%%%%%%%%%%%%%%%%%%%%%%%%%%%%%%%%%%

\begin{thebibliography}{999}
\bibitem{andrewsetal} Andrews, G.E.; Askey, R.; Roy, R. {Special functions}. 
In \emph{Encyclopedia of Mathematics and its Applications, \poro{71}};%MDPI: please confirm if the 71 can be deleted.
% AUTHOR: 71 should remain in the reference. URL https://mathscinet.ams.org/mathscinet/relay-station?mr=1688958
~Cambridge University Press: 
Cambridge, UK, 1999; pp. xvi+664.


%7
\bibitem{dlmf} Olver, F.W.J.; Daalhuis, A.B.O.; Lozier, D.W.; Schneider, B.I.; Boisvert, R.F.; Clark, C.W.; Miller, B.R.; Saunders, B.V.; Cohl, H.S.; McClain, M.A. NIST Digital Library of Mathematical Functions.  Available online: \href{https://dlmf.nist.gov/}{https://dlmf.nist.gov/}  \poro{(accessed on June 2023)}.%mdpi: please add the accessed date.
% AUTHOR: Done

%3
\bibitem{gas1989} Gasper, G. {$q$-extensions of Barnes', Cauchy's, and Euler's beta integrals}. 
In \emph{Topics in Mathematical Analysis}; Volume 11 of series Pure Maths; World Scientific Publishing: Teaneck, NJ, USA, 1989; pp. 294--314.

%4
\bibitem{gasrah}  Gasper, G.; Rahman, M. {Basic hypergeometric series}. 
\poro{With a foreword by Richard Askey.}
In \emph{Encyclopedia of Mathematics and its Applications}, 2nd ed.; 
Cambridge University Press: Cambridge, UK, 2004; Volume \poro{96.}%MDPI: We removed ``With a foreword by Richard Askey''. please confirm.
% AUTHOR: It should be on the reference. URL: https://mathscinet.ams.org/mathscinet/relay-station?mr=2128719

%9
\bibitem{sfa1} Sfaxi, R. On the Laguerre-Hahn intertwining operator
and application to connection formulae. {\em Acta Appl. Math.} \textbf{2011}, {\em 113}, 305--321.



%5
\bibitem{kolest} Koekoek, R.; Lesky, P.A.; Swarttouw, R.F. 
{Hypergeometric Orthogonal Polynomials and Their $q$-analogues}.
In \emph{Springer Monographs in Mathematics}; Springer: Berlin, Germany, 2010.

%6
\bibitem{leb} Lebedev, N.N.  {\em Special Functions and
Their Applications}; \poro{Prentice-Hall,  Englewood Cliffs, New Jersey }%mdpi: please add the location of the publisher.
% AUTHOR: Amended
1965. (In Russian)




%2
\bibitem{cohletall2020} Cohl, H.S.; Costas-Santos, R.S.
Multi-integral representations for associated Legendre and Ferrers functions. 
{\em Symmetry} {\bf 2020}, \mbox{\emph{12}, 22.}



%8
\bibitem{rebocho} Rebocho, M.N. Laguerre-Hahn orthogonal polynomials on the real line. 
{\em Random Matrices Theory Appl.} \textbf{2020}, {\em 9}, 33.


\end{thebibliography}
\end{document}